\documentclass[preprint,10pt]{elsarticle}

\usepackage{amsmath, amssymb}
\usepackage{graphicx}
\journal{Applied Mathematics and Computation}

\begin{document}

\begin{frontmatter}

\title{On the existence of eigenvalues of a boundary value problem with transmitting condition of the integral form for parabolic-hyperbolic equation}

\author{A.S.Berdyshev}

\address{Kazakh National Pedagogical University named after Abai, Almaty, Kazakhstan}
\ead{berdyshev@mail.ru}
\ead[url]{http://berdyshevas.narod.ru}

\author{E.T.Karimov}
\address{Institute of Mathematics, National University of Uzbekistan, Durmon yuli str.,29, 100125 Tashkent, Uzbekistan}
\ead{erkinjon@gmail.com}
\ead[url]{http://mathinst.uz/karimov/}

\begin{abstract}
In the paper, we investigate a local boundary value problem with transmitting condition of the integral form for mixed parabolic-hyperbolic equation with non-characteristic line of type changing. Two theorems on strong solvability and the existence of eigenvalues of the considered problem have been proved.
\end{abstract}


\end{frontmatter}

\section*{Introduction}

Main problems of the spectral theory of boundary value problems (BVPs) for mixed type equations one can divide as follows:

1)	characterization of the spectrum of boundary problems;

2)	construction of root (eigenfunctions and associated functions) functions;

3)	investigation of the completeness and basis property of root functions in various functional spaces.

Investigation of BVPs for mixed type equations becomes one of the main problems of the general theory of partial differential equations due to several applications of it in both in practice and theory. Nevertheless, despite the great attention to this problem by mathematicians, questions of the spectral theory of BVPs, in particular, for equations of mixed parabolic-hyperbolic type equations with integral transmitting conditions, remained hitherto unexplored.

In the work \cite{bcka}, analog of the generalized Tricomi problem with integral gluing conditions for mixed parabolic-hyperbolic equation was studied. Theorems on strong solvability and on the absence of eigenvalues were proved. In \cite{bcka} one can find historical information and notation on main scientific results on related field.

Omitting huge amount of work, we just note some of them, which are closely related to the present problem. One of the first investigations of BVPs with non-continuous transmitting conditions for parabolic-hyperbolic equations was work \cite{struch}. In \cite{ladst} authors investigated initial-boundary value problems for mixed type equations in multi-dimensional domains, which appear at studying problems on motion of conducting fluid in an electromagnetic field.

In the work \cite{ufl} the propagation of electrical oscillations in composite lines, when on the interval $0<x<l$ of the semi-infinite line losses are neglected, and the rest of the line is considered as a cable without leakage was reduced to solving a system of equations
\[\left.
\begin{array}{l}
  \displaystyle{ L\frac{\partial {{I}_{1}}}{\partial t}+\frac{\partial {{U}_{1}}}{\partial x}=0,\,\,\,\,{{C}_{1}}\frac{\partial {{U}_{1}}}{\partial t}+\frac{\partial {{I}_{1}}}{\partial x}=0,\,\,\,\,\,0<x<l\,} \\
\\
\displaystyle{   R{{I}_{2}}+\frac{\partial {{U}_{2}}}{\partial x}=0,\,\,\,\,{{C}_{2}}\frac{\partial {{U}_{2}}}{\partial t}+\frac{\partial {{I}_{2}}}{\partial x}=0,\,\,\,\,\,l<x<\infty}  \\
\end{array} \right\}
\]
with initial
\[
\left. U_1 \right|_{t=0}=0,\,\,\,\left. I_1 \right|_{t=0}=0,\,\,\,\left. U_2 \right|_{t=0}=0
\]
and boundary conditions
\[
\left. U_1 \right|_{x=0}=E(t),\,\,\,\underset{x\to \infty }{\mathop{\lim }}\,U_2=0,
\]
together with requirement of the continuity of the voltage and current
\[
\left. U_1 \right|_{x=l}=\left.U_2\right|_{x=l},\,\,\,\left. I_1 \right|_{x=l}=\left. I_2 \right|_{x=l}.
\]
Here $L,\,C_1$ are inductance and capacitance (per unit length) of the first part of the line; $R,\,C_2$ are resistance and capacitance of the second part.

Easy to ensure that if one exclude current from the system, the following parabolic-hyperbolic equation
\[
0=\left\{ \begin{array}{l}
   a_1^2 U_{xx}-U_{yy},\,\,\,\,0<x<l, \\
   \\
  a_2^2 U_{xx}-U_{y},\,\,\,\,l<x<+\infty  \\
\end{array} \right.
\]
can be deduced together with boundary conditions
\[
 U(x,0)=0,\,U_y(x,0)=0,\,\,\,0<x<l,\,\,U(x,0)=0,\,\,\,l<x<\infty,
\]
\[
U(0,y)=E(y),\,\,\underset{x\to +\infty }{\mathop{\lim }}\,U(x,y)=0.
\]
In this case transmitting condition will have form of
\[
U(l-0,y)=U(l+0,y),\,\,U_x(l+0,y)=\frac{R}{L}\int\limits_0^x{U_x(l-0,\eta )d\eta },
\]
where
\[
a_1^2=\frac{1}{LC_1},\,\,\,a_2^2=\frac{1}{RC_2}.
\]
Another example of application can be found in the work \cite{terl}.

Investigation of the unique solvability and spectral questions of some BVPs with integral transmitting conditions for parabolic-hyperbolic equations were done in works \cite{mois} -- \cite{kar}. Regarding the investigation of semilinear parabolic-equations we refer the readers to the works \cite{ferr} -- \cite{ashyr}.

In the present work, a new class of local BVPs with integral transmitting conditions for mixed parabolic-hyperbolic equations, which have eigenvalues, is found.

\section*{Formulation of the problem }
Let $\Omega \subset R^2$ be a finite simple-connected domain, bounded at $y>0$ by segments $AA_0, \, A_0B_0, \, B_0B \left(A=\left( 0,0 \right), A_0=\left( 0,1 \right), B_0=\left( 1,1 \right), B\left( 1,0 \right)\right)$, and at \\
$y<0$ by characteristics $AC:\,\,x+y=0$, $BC:\,\,x-y=1$ of the equation
\[
Lu=f\left( x,y \right),\eqno (1)
\]
where
\[
Lu\equiv \left\{ \begin{matrix}
   u_x-u_{yy},\,\,\,y>0,  \\
   u_{xx}-u_{yy},\,\,\,y<0.  \\
\end{matrix} \right.
\]
We consider the following variant of the Tricomi problem for parabolic-hyperbolic equation.

\textbf{Problem B.} Find a solution of the equation (1), satisfying boundary conditions
\[
u\left|_{AA_0\cup A_0B_0}=0, \right.\eqno (2)
\]
\[
u_x+u_y\left|_{BC}=0 \right.\eqno (3)
\]
and transmitting condition on the type-changing line
\[
u_x(x,+0)=u_x(x,-0),\,\,\,\,u_y(x,+0)=\alpha u_y(x,-0)-\beta \int\limits_0^x{u_y(t,-0)dt},\,\,0<x<1, \eqno (4)
\]
where $\alpha,\,\beta \in R$, such that $\alpha^2+\beta^2>0$.

Denote $\Omega_0=\Omega \cap \left\{ y>0 \right\},\,\Omega_1=\Omega \cap \left\{ y<0 \right\}$, let $\mathbb{W}$ be a set of functions belong to
\[
C\left( \bar{\Omega } \right)\cap C^{1,2}\left( \bar{\Omega }_0 \right)\cap C^{2,2}\left( \bar{\Omega}_1\right),
\]
satisfying equation (1) and condition (2)-(4).

Function $u\in L_2\left( \Omega  \right)$ we call as \emph{strong solution} of the problem, if there exists a sequence of functions $\left\{ u_n \right\}$, $u_n\in \mathbb{W}$, such that $\left\| u_n-u \right\|_{W_2^1(\Omega_0)}+\left\| u_n-u \right\|_{W_2^1(\Omega_1)}\to 0$, $\left\| Lu_n-f \right\|_0\to 0$ for $n\to \infty$. Here and further, by $\left\| \cdot  \right\|_l$ we denote norm of the Sobolev space $W_2^l\left( \Omega  \right)$, where $W_2^0\left( \Omega  \right)=L_2\left( \Omega  \right)$.

\section*{Main result}

\textbf{Theorem 1.} For any function $f\in L_2\left( \Omega  \right)$ there exists unique strong solution of the problem B. This solution belongs to the class of functions $W_2^1\left( \Omega  \right)\cap W_2^{1,2}\left( \Omega_0 \right)\cap C\left( \bar{\Omega} \right)$, satisfies the following inequality
\[
\left\| u \right\|_{W_2^1(\Omega_0)}+\left\| u \right\|_{W_2^1(\Omega_1)}\le c\left\| f \right\|_0\eqno (5)
\]
and represented as
\[
u\left( x,y \right)=\iint\limits_{\Omega }{K\left( x,y; x_1, y_1\right)f\left(x_1, y_1\right)dx_1dy_1},\eqno (6)
\]
where $K\in L_2\left( \Omega \times \Omega  \right)$.

\textbf{Proof}: According to the unique solvability of the first boundary problem for the heat equation with conditions (2) and $u\left( x,0 \right)=\tau \left( x \right)$, $\tau \left( 0 \right)=0$, and the Darboux problem for the wave equation with conditions (3), $u\left( x,0 \right)=\tau \left( x \right)$, $\tau \left( 0 \right)=0$, solution of the equation (1) can be represented as follows
\[
u\left( x,y \right)=\left\{
\begin{array}{l}
   \int\limits_0^x dx_1\int\limits_0^1 G\left( x-x_1,y,y_1\right)f\left( x_1,y_1\right)dy_1 +\\
   +\int\limits_0^x G_{y_1}\left( x-x_1,y,0 \right)\tau \left(x_1\right)dx_1,\,\,\,y>0, \\
   \\
   \int\limits_{\xi}^{\eta } d\xi_1\int\limits_{\eta}^1 f_1\left( \xi_1,\eta_1\right)d\eta_1+\tau \left( \eta  \right),\,\,\,y<0,\\
\end{array} \right.\eqno (7)
\]
where $f_1\left( \xi ,\eta  \right)=\frac{1}{4}f\left( \frac{\xi +\eta }{2},\frac{\xi -\eta }{2} \right)$, $\xi =x+y$, $\eta =x-y$, and $G\left( x-x_1,y,y_1\right)$ is Green's function of the first boundary problem for heat equation in a rectangle $AA_0B_0B$, which has a form:
\[
G\left( x,y,y_1\right)=\frac{1}{2\sqrt{\pi x}}\sum\limits_{n=-\infty }^{+\infty }{\left[ e^{ -\frac{\left( y-y_1+2n \right)^2}{4x}}-e^{ -\frac{\left( y+y_1+2n \right)^2}{4x}} \right]}.\eqno (8)
\]

Considering (8), calculating derivative $\frac{\partial u}{\partial y}$, and passing to the limit as $y$ tends to zero in $\Omega_0$, we obtain first functional relation between functions $\tau \left( x \right)$ and $\nu_1\left( x \right)=\frac{\partial u}{\partial y}\left( x,+0 \right)$ given as
\[
\nu_1\left( x \right)=-\int\limits_0^x k\left( x-t \right)\tau'\left(t\right)dt+ \Phi_0\left( x \right),\eqno (9)
\]
where
\[
k\left( x \right)=\frac{1}{\sqrt{\pi x}}\sum\limits_{n=-\infty }^{+\infty }{e^{-\frac{n^2}{x}}},\eqno (10)
\]
\[
\Phi_0\left( x \right)=\int\limits_0^x dx_1\int\limits_0^1 G_0\left( x-x_1,y_1\right)f\left(x_1,y_1\right)dy_1,\eqno (11)
\]
\[
G_0\left( x,y_1\right)\equiv G_y\left( x,y_1,0 \right)=\frac{1}{2\sqrt{\pi } x^{3/2}}\sum\limits_{n=-\infty }^{+\infty }{\left(y_1+2n \right) e^{-\frac{\left( y_1+2n \right)}{4x}}}.\eqno (12)
\]

Similarly, we find another integral-differential relation between functions $\tau \left( x \right)$ and $\nu_2\left( x \right)=\frac{\partial u}{\partial y}\left( x,-0 \right)$ on $AB$, reduced from the domain $\Omega_1$. It has the form
\[
\nu_2\left( x \right)=-\tau'\left( x \right)-2\int\limits_x^1 f_1\left( x,t \right)dt,\,\, 0<x<1.\eqno (13)
\]
Let $\alpha \ne 0$. From (9) and (13) based on transmitting conditions (4), we deduce integral equation with respect to $\tau'(x)$:
\[
\tau'(x)+\int\limits_0^x k_1(x-t)\tau'(t)dt=F(x).\eqno (14)
\]
Here
\[
k_1(x-t)=\frac{1}{\alpha }\left[ k(x-t)+\beta \right],\eqno (15)
\]
\[
F(x)=-\frac{1}{\alpha }\Phi_0(x)-2\int\limits_x^1 f_1(x,\eta_1)d\eta_1+\frac{2\beta }{\alpha }\int\limits_0^x dt\int\limits_t^1 f_1(x,\eta _1)d\eta _1.\eqno (16)
\]
Thus, the problem is equivalently reduced to the second kind Volterra integral equation (14). Since by (10), the kernel $k\left( x \right)$ can be represented as
\[
k\left( x \right)=\frac{1}{\sqrt{\pi x}}+\tilde{k}\left( x \right),
\]
where $\tilde{k}\left( x \right)\in C^{\infty }\left[ 0;1 \right]$, then from (15) it follows that $k_1\left( x \right)$ has weak singularity. Therefore, there exists the unique solution of (14) and it has a form
\[
\tau'\left( x \right)=F\left( x \right)+\int\limits_0^x \Gamma \left( x-t \right)F\left( t \right)dt,\eqno	(17)
\]
where $\Gamma \left( x \right)$  is resolvent kernel of (14), which defined as
\[
\Gamma \left( x \right)=\sum\limits_{j=1}^{\infty }{k_j\left( x \right)},
\]
\[
k_{j+1}\left( x \right)=\int\limits_0^x{k_1\left( x-t \right) k_j\left( t \right)dt},\, j\in N.
\]
Considering $\tau \left( 0 \right)=0$ after some evaluations, from (17) we get
\[
\begin{array}{l}
   \tau \left( x \right)=-\frac{1}{\alpha }\int\limits_{0}^{x}{d{{x}_{1}}}\int\limits_{0}^{1}{{{G}_{1}}\left( x-{{x}_{1}},{{y}_{1}} \right)f}\left( {{x}_{1}},{{y}_{1}} \right)d{{y}_{1}}- \\
  \\
  -2\int\limits_{0}^{x}{d{{\xi }_{1}}\int\limits_{{{\xi }_{1}}}^{1}{{{\Gamma }_{1}}\left( x,{{\xi }_{1}} \right){{f}_{1}}\left( {{\xi }_{1}},{{\eta }_{1}} \right)d{{\eta }_{1}}+\frac{2\beta }{\alpha }\int\limits_{0}^{x}{d{{\xi }_{1}}}\int\limits_{{{\xi }_{1}}}^{1}{{{f}_{1}}}\left( {{\xi }_{1}},{{\eta }_{1}} \right)d{{\eta }_{1}}\int\limits_{{{\xi }_{1}}}^{x}{{{\Gamma }_{1}}}\left( x,t \right)dt}} \\
\end{array}\eqno  (18)
\]
where
\[
\Gamma_1\left( x,t \right)=1+\int\limits_t^x \Gamma \left( z-t \right)dz, \eqno	(19)
\]
\[
G_1\left( x-x_1,y_1 \right)=\int\limits_{x_1}^x {G_0\left( t-x_1,y_1\right)\Gamma_1\left( x,t \right)dt}.
\]
Substituting (18) into (7), we deduce the formula (6), where
\[
\begin{array}{l}
   K\left( x,y;x_1,y_1\right)=\theta \left( y \right)\left\{ \theta \left( y_1\right)\theta \left( x-x_1\right)G_2\left( x-x_1,y,y_1\right)- \right.\\
  \\
   \left. -\theta \left( -y_1\right)\left[ \theta \left( x-\xi_1\right)G_3\left( x-\xi_1,y \right)+\frac{\beta }{\alpha }\theta \left( x-\xi_1\right)G_4\left( x-\xi_1,y \right) \right] \right\}+\\
  \\
   +\theta \left( -y \right)\left\{ -\frac{1}{\alpha }\theta \left( y_1\right)\theta \left( \eta -x_1 \right)G_1\left( \eta -x_1,y_1\right)+\theta \left( -y_1\right)\left[ \frac{1}{2}\theta \left( \eta -\xi_1\right)\times\right.\right.\\
  \\
   \left.\left.\times\theta \left( \xi_1-\xi  \right)\theta \left(\eta_1-\eta  \right)-\theta \left( \eta -\xi_1\right)\left[ \Gamma_1\left( \eta ,\xi_1\right)+\frac{\beta }{\alpha }\int\limits_{\xi_1}^{\eta }{\Gamma \left( \eta -t \right)dt} \right] \right]\right\}, \\
\end{array}\eqno (20)
\]
\[
G_2\left( x-x_1,y \right)=G\left( x-x_1,y,y_1 \right)-\frac{1}{\alpha }\int\limits_{x_1}^x G_{y_1}\left( x-t,y,0 \right) G_1\left( t-x_1,y_1 \right)dt,
\]
\[
G_3\left( x-\xi_1,y \right)=\int\limits_{\xi_1}^x G_{y_1}\left( x-t,y,0 \right)\Gamma_1\left( t, \xi_1\right)dt,
\]
\[
G_4\left( x-\xi_1,y \right)=\int\limits_{\xi_1}^x dt\int\limits_{\xi_1}^x G_{y_1}\left( x-t,y,0 \right) \Gamma_1\left( t,z \right)dz,
\]
\[
\theta \left( y \right)=\left\{ \begin{array}{l}
   1,\,\,\,\,\,if\,\,\,\,y>0, \\
  0,\,\,\,if\,\,\,\,y<0. \\
\end{array} \right.
\]

Similarly as in [1] (see proof of the Lemma 1), one can prove that
\[K\left( x,y;x_1,y_1\right)\in L_2\left( \Omega \times \Omega  \right).\]
Considering (11) by virtue of direct calculations from (16) we can state that the estimate
\[
\left\| F\left( x \right) \right\|_{L_2\left( 0,1 \right)}\le C\left\| f \right\|_0
\]
is valid. Therefore, from (14) we have
\[
\left\| \tau'\left( x \right) \right\|_{L_2\left( 0,1 \right)}\le C\left\| F\left( x \right) \right\|_{L_2\left( 0,1 \right)}\le C\left\| f \right\|_0.
\]
Based on this and properties of the solution of the first boundary problem for heat equation, it follows that solution of the problem B belongs to the class of functions $W_2^1\left( \Omega  \right)\cap W_2^{1,2}\left(\Omega_0\right)\cap C\left( \bar{\Omega} \right)$ and satisfies inequality (5).

Now we show that found solution will be strong. Since $C_0^1\left(\bar{\Omega} \right)$ is dense in $L_2\left( \Omega  \right)$, then for any function $f\in L_2\left( \Omega  \right)$ there exists a sequence of functions $f_n\in C_0^1\left( \bar{\Omega } \right)$ such that $\left\|f_n-f \right\|\to 0$, $n\to 0$. Here $C_0^1\left(\bar{\Omega} \right)$ is a set of differentiable functions in $\Omega$, which are equal to zero in neighborhood of $\partial \Omega$ ($\partial \Omega$ is a boundary of the domain $\Omega$). Denote that $u_n=\mathfrak{L}^{-1} f_n$ ,where  $\mathfrak{L}^{-1}$ is inverse of the operator $\mathfrak{L}$ of the problem �.

Easy to conclude that $F_n\left( x \right)\in C^1\left[ 0;1 \right]$ at $f_n\in C_0^1\left( \bar{\Omega} \right)$. Here by $F_n(x)$ we denote representation similar to (16), where $f(x,y)$ should be replaced with $f_n(x,y)$. Hence, equation (14) we can consider as the second kind Volterra integral equation in the space $C^1\left[ 0;1 \right]$. Consequently, $\tau'_n\left( x \right)=u_{nx}\left( x,0 \right)\in C^1\left[ 0;1 \right]$. Based on properties of the first boundary problem for heat equation and the Darboux problem for wave equation, considering representation (6), we get that $u_n\in \mathbb{W}$ for all $f_n\in C_0^1\left( \bar{\Omega } \right)$.

By virtue of the inequality we obtain
\[
\left\|u_n-u \right\|_{W_2^1(\Omega_0)}+\left\| u_n-u \right\|_{W_2^1(\Omega_1)}\le c\left\|f_n-f \right\|_0\to 0.
\]
Therefore, $\left\{ u_n\right\}$ sequence satisfies all requirements of the definition of the strong solution. Now we can state that the problem B is strongly solvable for any $f$, and strong solution belongs to the class of functions $W_2^1\left( \Omega  \right)\cap W_2^{1,2}\left( \Omega_0 \right)\cap C\left( \bar{\Omega } \right)$.

Theorem 1 is proved.

From the Theorem 1 we can conclude that operator $\mathfrak{L}$ of the problem B is invertible, and inverse operator $\mathfrak{L}^{-1}$ is Hilbert-Schmidt operator. There is a natural question on the existence of eigenvalues of the operator $\mathfrak{L}^{-1}$, consequently, of the problem B as well.

\textbf{Theorem 2.} Let $\alpha >0,\,\beta >0$. There exists $\lambda \in \mathbb{C}$ such that the equation
\[
Lu=\lambda u
\]
has non-trivial solution $u\in \mathbb{W}$.

\textbf{Proof}: We denote by $\mathfrak{L}$ a closure in $L_2\left( \Omega  \right)$ of the differential operator given in $\mathbb{W}$ by equality (1). From the theorem 1 follows that $\mathfrak{L}$ is invertible and $\mathfrak{L}^{-1}$ defined by (6) is the Hilbert-Schmidt operator. Then operator $\mathfrak{L}^{-2}\equiv \left( \mathfrak{L}^{-1} \right)^2$ is nuclear in $L_2\left( \Omega  \right)$. Therefore, we apply the result of the Lidskii \cite{lids} on a coincidences of the matrix and spectral traces to the operator $\mathfrak{L}^{-2}$.

\textbf{Lemma 1.} \cite{lids} If operator $T$ is nuclear in a Hilbert space $H$, then for any orthonormal basis $\phi_i\left( i=1,2,... \right)$ in $H$, the equality
\[
SpT\equiv \sum\limits_{k=1}^{\infty }{\left( T\phi_k, \phi_k\right)}=\sum\limits_{k=1}^{\infty }{\lambda_k\left( T \right)}
\eqno (21)
\]
holds true. Here $\lambda_k$ are eigenvalues of the operator $T$.

Known that if operator $T$ is nuclear in $L_2\left( \Omega  \right)$, represented as multiplication $T=PR$ of two Hilbert-Schmidt operators
\[
\left( Pf \right)\left( z \right)=\int\limits_{\Omega }{P\left( z,z_1 \right)f\left( z_1 \right)dz_1},
\]
\[\left( Rf \right)\left( z \right)=\int\limits_{\Omega }{R\left( z,z_1\right)f\left( z_1 \right)dz_1},
\]
then Gaal's formula for calculating traces \cite{bris}
\[
SpT=\int\limits_{\Omega }{\left[ \int\limits_{\Omega }{P\left( z,z_1 \right)R\left( z_1,z \right)dz_1} \right]}dz\eqno (22)
\]
is true.

From (21) and (22) we deduce

\[
Sp \mathfrak{L}^{-2}=\iint\limits_{\Omega }{dxdy\iint\limits_{\Omega }{K\left( x,y;x_1,y_1\right)K\left( x_1,y_1;x,y \right)dx_1dy_1}}.
\]
From (20) it follows that
\[
\begin{array}{l}
K\left( x,y;x_1,y_1 \right)K\left(x_1,y_1;x,y \right)=\theta \left( y \right)\theta \left( -y_1 \right)\theta \left( x-\xi_1\right)\theta \left( \eta_1-x \right)\times\\
\\
\times\left[ \frac{1}{\alpha }G_1(\eta_1-x,y) G_3(x-\xi_1,y)+
\frac{\beta }{\alpha }G_1\left( \eta_1-x,y \right)G_4\left( x-\xi_1,y \right) \right]+\\
\\
+\theta \left( -y \right)\theta \left(y_1 \right)\theta \left( \eta -x_1\right)\theta \left( x_1-\xi  \right)\left[ \frac{1}{\alpha }G_1(\eta -x_1,y_1) \right.G_3(x_1-\xi ,y_1)+\\
\\
+\left. \frac{\beta }{\alpha^2}G_1\left( \eta -x_1,y_1 \right)G_4\left(x_1-\xi ,y_1 \right) \right]+\theta \left( -y \right)\theta \left( -y_1 \right)\left\{ -\frac{1}{2}\theta \left( \eta -\xi_1\right)\theta \left(\xi_1-\xi  \right)\right.\times\\
\end{array}
\]

\[
\begin{array}{l}
\times \theta \left(\eta_1-\eta  \right)\theta \left(\eta _1-\xi  \right)\left[ \Gamma_1\left(\eta_1,\xi  \right)-\frac{\beta }{\alpha }\int\limits_{\xi_1}^{\eta }\Gamma\left( \eta -t \right)dt \right]d\eta_1-\frac{1}{2}\theta \left( \eta -\xi_1\right)\times\\
\\
\times \theta \left(\eta_1-\xi  \right)\theta \left( \xi -\xi_1\right)\theta \left( \eta -\eta_1\right)\left[\Gamma_1\left(\eta_1,\xi  \right)-\frac{\beta }{\alpha }\int\limits_{\xi_1}^{\eta }\Gamma\left( \eta -t \right)dt \right]+\\
\\
+\theta \left( \eta -\xi_1 \right)\theta \left( \eta_1-\xi  \right)\left[\Gamma_1\left(\eta_1,\xi  \right)\Gamma_1\left( \eta ,\xi_1 \right)+\frac{\beta }{\alpha }\Gamma_1\left(\eta_1,\xi  \right)\int\limits_{\xi_1}^{\eta }\Gamma\left( \eta -t \right)dt \right.+\\
\\
+\left. \left. \frac{\beta }{\alpha }\Gamma_1\left( \eta ,\xi_1\right)\int\limits_{\xi }^{\eta_1}\Gamma\left(\eta_1-t \right)dt+\left( \frac{\beta }{\alpha } \right)^2\int\limits_{\xi }^{\eta_1}\Gamma\left(\eta_1-t \right)dt\int\limits_{\xi_1}^{\eta }\Gamma\left( \eta -t \right)dt \right] \right\}.\\
\end{array}
\]
Therefore,
\[
\begin{array}{l}
Sp \mathfrak{L}^{-2}=\frac{2}{\alpha }\iint\limits_{\Omega_1}dxdy\iint\limits_{\Omega_2}\theta \left( x-\xi_1\right)\theta \left(\eta_1-x \right)\left[G_1(\eta_1-x,y)G_4(x-\xi_1,y)+\right.\\
\\
\left.+\frac{\beta }{\alpha }G_1\left(\eta_1-x,y \right)G_3\left( x-\xi_1,y \right) \right]dx_1dy_1-\iint\limits_{\Omega_2}dxdy\iint\limits_{\Omega_2}\theta \left(\eta -\xi_1\right)\times\\
\\
\times\theta \left(\xi_1-\xi  \right)\theta \left(\eta_1-\eta  \right)\theta \left(\eta_1-\xi  \right)\left[\Gamma_1\left(\eta_1,\xi  \right)+\frac{\beta }{\alpha }\int\limits_{\xi_1}^{\eta }\Gamma\left( \eta -t \right)dt \right]dx_1dy_1+\\
\\
  +\iint\limits_{\Omega_2}dxdy\iint\limits_{\Omega_2}\theta \left( \eta -\xi_1 \right)\theta \left(\eta_1-\xi  \right)\left[ \Gamma_1\left( \eta_1,\xi  \right)\Gamma_1\left( \eta ,\xi_1 \right)+\frac{\beta }{\alpha }\Gamma_1\left( \eta_1,\xi  \right)\right.\times\\
\\
  \times\int\limits_{\xi_1}^{\eta }\Gamma\left(\eta -t \right)dt + \frac{\beta }{\alpha }\Gamma_1\left( \eta ,\xi_1 \right)\int\limits_{\xi }^{\eta_1}\Gamma\left(\eta_1-t \right)dt+\\
\\
 +\left. \left( \frac{\beta }{\alpha } \right)^2\int\limits_{\xi }^{\eta_1}\Gamma\left(\eta_1-t \right)dt\int\limits_{\xi_1}^{\eta }\Gamma\left( \eta -t \right)dt \right]dx_1dy_1=\sum\limits_{k=1}^3 I_k. \\
\end{array}\eqno (23)
\]
Let us show that $\sum\limits_{k=1}^3 I_k>0$. In fact
\[
\begin{array}{l}
  I_3+I_2=\frac{1}{4}\int\limits_0^1d\eta \int\limits_0^{\eta }d\xi \int\limits_0^1d\xi_1\int\limits_{\xi_1}^1\theta \left( \eta -\xi_1 \right)\theta \left(\eta_1-\xi  \right)\left[\Gamma_1\left( \eta ,\xi_1\right)\Gamma_1\left(\eta_1,\xi  \right)\right.+\\
\\
  +\frac{\beta }{\alpha }\Gamma _1\left(\eta_1,\xi  \right)\int\limits_{\xi_1}^{\eta }\Gamma\left( \eta -t \right)dt
 +\frac{\beta }{\alpha }\Gamma_1\left( \eta ,\xi_1 \right)\int\limits_{\xi }^{\eta_1}\Gamma\left(\eta_1-t \right)dt+\\
\\
 +\left.\left( \frac{\beta }{\alpha } \right)^2\int\limits_{\xi }^{\eta_1}\Gamma\left(\eta_1-t \right)dt\int\limits_{\xi_1}^{\eta }{\Gamma }\left( \eta -t \right)dt \right]- \frac{1}{4}\int\limits_0^1 d\eta \int\limits_0^{\eta }d\xi \int\limits_0^1 d\xi_1\int\limits_{\xi_1}^1 \theta \left( \eta -\xi_1\right)\times\\
\\
 \times\theta \left(\xi_1-\xi  \right)\theta \left(\eta_1-\eta  \right)\theta \left(\eta_1-\xi  \right)\left[\Gamma_1\left(\eta_1,\xi  \right)+\frac{\beta }{\alpha }\int\limits_{\xi_1}^{\eta }\Gamma\left( \eta -t \right)dt \right]d\eta_1=\\
\\
=\frac{1}{4}\int\limits_0^1d\eta \int\limits_0^{\eta }d\xi \int\limits_0^{\eta }d\xi_1\int\limits_{\xi_1}^1\theta \left(\eta_1-\xi  \right)\left[\Gamma_1\left(\eta_1,\xi  \right)+\frac{\beta }{\alpha }\int\limits_{\xi_1}^{\eta }\Gamma\left( \eta -t \right)dt \right]\times \\
\\
\times \left\{ \Gamma_1\left( \eta ,\xi_1 \right)+\frac{\beta }{\alpha }\int\limits_{\xi }^{\eta }\Gamma\left(\eta_1-t \right)dt-\theta \left( \xi_1-\xi  \right)\theta \left( \eta_1-\eta  \right) \right\}d\eta_1>0,\\
\end{array}\eqno (24)
\]
Since, due to (8), (15), (19)  $\Gamma_1\left( \eta ,\xi_1 \right)\ge 1$ and $\alpha >0,\,\,\beta >0$, then
\[
\Gamma_1\left( \eta ,\xi_1\right)+\frac{\beta }{\alpha }\int\limits_{\xi_1}^{\eta }\Gamma\left( \eta -t \right)dt-\theta \left(\xi_1-\xi  \right)\theta \left(\eta_1-\eta  \right)>0.
\]
Now consider $I_1$. We have
\[
\begin{array}{l}
I_1=\frac{2}{\alpha }\int\limits_0^1 dx\int\limits_0^1 dy\int\limits_0^1 d\xi_1\int\limits_{\xi_1}^1 \theta \left( x-\xi_1 \right)\theta \left( \eta_1-x \right)G_1(\eta_1-x,y)\left[G_4(x-\xi_1,y)+\right.\\
\\
\left.+\frac{\beta }{\alpha }G_3\left( x-\xi_1,y \right) \right]d\eta _1=\frac{2}{\alpha }\int\limits_0^1 dx\int\limits_0^1 dy\int\limits_0^x d\xi_1\int\limits_x^1 G_1\left(\eta_1-x,y \right)\left[G_4(x-\xi_1,y)+\right.\\
\\
\left.+\frac{\beta }{\alpha }G_3\left( x-\xi_1,y \right) \right]d\eta_1=\frac{2}{\alpha }\int\limits_0^1 dx\int\limits_0^1 dy\int\limits_0^x d\xi_1\int\limits_0^{1-x}G_1\left( x-\xi_1,y \right)\times\\
\\
\left[ G_4(\eta _2,y)+\frac{\beta }{\alpha }G_3\left(\eta_2,y \right) \right]d\eta_2=\frac{2}{\alpha }\int\limits_0^1 dx\int\limits_0^1 dy\left( \int\limits_0^x G_1\left(\xi_2,y \right)d\xi_2\right)\times\\
\\
\left( \int\limits_0^{1-x}\left[G_4(\eta_2,y)+\frac{\beta }{\alpha }G_3\left(\eta_2,y \right) \right]d\eta_2 \right).\\
\end{array}\eqno (25)
\]
Function $G_1$ we represent as
\[
G_1\left( \xi ,y \right)=\int\limits_0^{\xi }G_0\left( t,y \right)dt+\int\limits_0^{\xi }G_0\left( t,y \right)dt\int\limits_0^{\xi -t}\Gamma \left( \tau  \right)d\tau.\eqno (26)
\]
Taking (12) into account, investigate first item. For this, we use the following transformations:
\[
\begin{array}{l}
\int\limits_0^{\xi }G_0\left( t,y \right)dt=\frac{1}{2\sqrt{\pi }}\sum\limits_{n=-\infty }^{\infty }{\int\limits_0^{\xi }\frac{y+2n}{t^{3/2}}e^{-\frac{\left( y+2n \right)^2}{4t}}dt}=\frac{2}{\sqrt{\pi }}\sum\limits_{n=-\infty }^{\infty }{\int\limits_{\frac{y+2n}{2\sqrt{\xi }}}^{\pm \infty }e^{-t^2}dt=}\\
\\
=-\frac{2}{\sqrt{\pi }}\sum\limits_{n=-\infty }^{-1}{\int\limits_{-\infty }^{\frac{y+2n}{2\sqrt{\xi }}}e^{-t^2}dt}+\frac{2}{\sqrt{\pi }}\sum\limits_{n=0}^{\infty }{\int\limits_{\frac{y+2n}{2\sqrt{\xi }}}^{+\infty }e^{-t^2}dt}=\\
\\
=\frac{2}{\sqrt{\pi }}\left[ \sum\limits_{n=0}^{\infty }{\int\limits_{\frac{y+2n}{2\sqrt{\xi }}}^{\infty }e^{-t^2}dt}- \sum\limits_{n=-1}^{\infty }{\int\limits_{-\frac{y+2n}{2\sqrt{\xi }}}^{\infty }e^{-t^2}dt} \right]=\\
\\
=\frac{2}{\sqrt{\pi }}\left[ \sum\limits_{n=0}^{\infty }{\int\limits_{\frac{y+2n}{2\sqrt{\xi }}}^{\infty }e^{-t^2}dt}-
\sum\limits_{n=0}^{\infty }{\int\limits_{\frac{2n-y+2}{2\sqrt{\xi }}}^{\infty }e^{-t^2}dt} \right]=\frac{2}{\sqrt{\pi }}\sum\limits_{n=0}^{\infty }{\int\limits_{\frac{2+yn}{2\sqrt{\xi }}}^{\frac{2n+2-y}{2\sqrt{\xi }}}e^{-t^2}dt}.\\
\end{array}
\]
From here we get
\[
\int\limits_{0}^{\xi }G_0\left( t,y \right)dt\ge 0.\eqno (27)
\]
The equality in (27) will be true only when $y=1$, i.e.
\[
\int\limits_{0}^{\xi }G_0\left( t,y \right)dt\not{\equiv }0.
\]
Then considering $\Gamma \left( \tau  \right)>0$, with respect to the second item of (26) we have
\[
\int\limits_0^{\xi }G_0\left( t,y \right)dt\int\limits_0^{\xi -t}\Gamma \left( \tau  \right)d\tau =\int\limits_0^{\xi }\Gamma \left( \tau  \right)d\tau \int\limits_0^{\xi -\tau }G_0\left( t,y \right)dt\ge 0\,\, \left( \not{\equiv }0 \right).
\]
Similarly, we can prove that the second item of (25) is as well positive. Hence, from (25) we can state that $I_1>0$. From (23)-(25) it follows that $Sp \mathfrak{L}^{-2}>0$.

Then by virtue of (20), we have
\[
\sum\limits_{k=1}^{\infty }\lambda_k\left(\mathfrak{L}^{-2}\right)\equiv \sum\limits_{k=1}^{\infty }\lambda_k^2\left(\mathfrak{L}^{-1} \right)>0,
\]
where $\lambda_k\left(\mathfrak{L}^{-2} \right)$ are eigenvalues of $\mathfrak{L}^{-2}$. It means that $\sum\limits_{k=1}^{\infty }{\frac{1}{\lambda_k^2}}>0$, where $\lambda_k$ are eigenvalues of the problem (1)-(3). From here, the existence of eigenvalues of the problem B follows.

Theorem 2 is proved.


\begin{thebibliography}{10}

\bibitem{bcka}	Berdyshev A.S., Cabada A., Karimov E.T. and Akhtaeva N.S. On the Volterra property of a boundary problem with integral gluing condition for a mixed parabolic-hyperbolic equation. Boundary Value Problems \textbf{2013}:94.

\bibitem{struch} Struchina G.M. A problem of pairing two equations. Injenerno-fizicheskiy jurnal \textbf{11} (4) (1961),99-104 .

\bibitem{ladst} Ladizhenskaya O.A., Stupyalis L. On mixed type equations. Vestnik LGU, seriya matematika, mexanika i astronomia \textbf{19} (4) (1965), 38-46.

\bibitem{ufl} Uflyand Y.S. On the question of the distribution of fluctuations in composite electrical lines. Injenerno-fizicheskiy jurnal \textbf{7} (1) (1964), 89-92.

\bibitem{terl} Terlyga O., Bellout H., Bloom F. A hyperbolic-parabolic system arising in pulse combustion: existence of solutions for the linearized problem. Electronic Journal of differential Equations, \textbf{2013} (46) (2013), 1-42.

\bibitem{mois} Moiseev E.I., Kapustin N.Y.: On spectral problems with a spectral parameter in the boundary conditions.Differ. Uravn. (Minsk) \textbf{33} (1) (1997), 115-119  (Russian), translation in Differ. Equ. \textbf{33} (1) (1997), 116-120.

\bibitem{sadt} Sadybekov M.A., Tojzhanova G.D. Spectral properties of a class of parabolic-hyperbolic equations. Differ. Uravn. (Minsk) \textbf{28} (1) (1992), 176-179.

\bibitem{ak} Akhtaeva N.S. Karimov E.T. On a boundary problem with gluing conditions of integral form for mixed parabolic-hyperbolic equation with non-characteristic line of type changing. Vestnik KazNU, seriya matem., mexan., inform \textbf{77}(2)  (2013), 64-77.

\bibitem{berd}	Berdyshev A.S. On the existence of eigenvalues of a boundary problem for parabolic-hyperbolic equation of the third order. Uzb. Math. Journal \textbf{2} (1998), 19-25.

\bibitem{kar} 	Karimov E.T. Non-local problems with special gluing condition for the parabolic-hyperbolic type equation with complex spectral parameter. PanAmerican Mathematical Journal \textbf{17}(2) (2007), 11-20.

\bibitem{ferr} Ferreira J. On weak solutions of semilinear hyperbolic-parabolic equations. Internat. J. Math. @ Math. Sci. \textbf{19}(4) (1996), 751-758.

\bibitem{ashyr}	Ashyralyev A., Yurtsever A. On a nonlocal boundary value problem for semilinear hyperbolic-parabolic equations. Nonlinear Analysis \textbf{47} (2001), 3585-3592.

\bibitem{lids} Lidskii V.B. Non-adjoin operators with traces. Doklady AN SSSR \textbf{125}(3) (1959), 485-488.

\bibitem{bris}	Brislawn C. Kernels of trace class operators. Proc. Amer. Math. Soc. \textbf{104}(4) (1988), 1181-1190.

\end{thebibliography}
\end{document}